\documentclass[12pt]{article}
\usepackage{fullpage}
\usepackage{amsmath}
\usepackage{amssymb,amsfonts,amsthm}
\usepackage{verbatim}

\newcommand{\KSS}{\ensuremath{\mbox{\sf KSS}}}

\newcommand{\STS}{\ensuremath{\mbox{\sf STS}}}

\theoremstyle{definition}
\newtheorem{construction}{Construction}

\theoremstyle{plain}

\newtheorem{theorem}[construction]{Theorem}

\theoremstyle{remark}

\title{Technical report\\
Existence of Kirkman signal sets\\
on $v=1,3\pmod{6}$ points, $14\leq v \leq 3000$.
}

\author
{
Melissa S. Keranen \quad Donald L. Kreher \\[2ex]
\it Department of Mathematical Sciences\\
\it Michigan Technological University\\
\it Houghton, MI 49931-0402, USA
}

\begin{document}
\maketitle

\renewcommand{\theconstruction}{\Alph{construction}}

A partial Steiner triple 
system whose triples can be partitioned into $s$ partial 
parallel classes, each of size $m$, is a {\em signal set}, 
denoted $\SS(v,s,m)$.  A {\em Kirkman signal set} 
$\KSS(v,m)$ is an $\SS(v,s,m)$ with 
$s=\lfloor\mu(v)/m\rfloor$. When $v \equiv 1$ or $3 
\pmod{6}$, then $\mu(v)=b$, so the decomposition of an 
$\STS(v)$ into partial parallel classes of size $m$ is 
equivalent to a $\KSS(v,m)$.  
In 2012, Colbourn, Horsley, and Wang 
\cite{CHW} established the following:

\begin{theorem}\emph{\cite{CHW}} \label{CHW0}
Let $u\equiv 2 \pmod{6}$ and $v\equiv 1,3 \pmod{6}$ be 
integers 
such  that $u \geq 88$ and $v \geq \textsc{max}(9u2+ 9u, 
99u)$. 
Then there exists an
 $\STS(v)$  whose triples may be
 ordered 
such that any  
\[
M(u,v)=
\left\lfloor \frac{1}{3}\bigg(\frac{u-2}{u}\bigg)
\bigg( \frac{2(v-6u+1)^2}{2v+6u^2-9u+2} \bigg) \right\rfloor
\]
consecutive triples  are pairwise vertex disjoint.
\end{theorem}
\noindent{}It follows that for all such $v,u$, that a $\KSS(v,m)$
will exist whenever $m \leq M(u,v)$.  
In the same paper, they showed 
\begin{theorem}\emph{\cite{CHW}}
If $14 \leq v \leq 32$ 
and $m \leq \lfloor v/3 \rfloor$, there exists a 
$\KSS(v,m)$.  
\end{theorem}
\noindent{}Hodaj, Keranen,  Kreher and Tollefson establish the following
new constructions in ~\cite{HKKT}.
\begin{theorem}\emph{\cite{HKKT}}
\label{m=4}
Suppose $v \equiv 1, 3 \pmod{6}$, $v \geq 19$, and $v \not 
= 21,27$. If $4|b$, then there exists a $\KSS(v,4)$.
\end{theorem}
\begin{theorem}\emph{\cite{HKKT}}
\label{C 1}
If $v \equiv 3 \pmod{6}$, $v>9$, then there exists a 
$\KSS(4v-3, v-1)$.
\end{theorem}
\begin{theorem}\emph{\cite{DR}(see also \cite{HKKT})}
\label{CII}
If $v \equiv 1 \pmod{6}$, then there exists a $KSS(v, 
\frac{v-1}{6})$.
\end{theorem}
\begin{theorem}\emph{\cite{HKKT}}
\label{g^2}
If $g \equiv 3 \pmod{6}$, then there exists a 
$\KSS(g^2, g(g-1)/3)$
\end{theorem}
\begin{theorem}\emph{\cite{HKKT}}
\label{gr}
Let $v=g \cdot r$ be a positive integer such that $g \equiv 
r \equiv 3 \pmod{6}$.  If $d$ is chosen so that $d | 
\frac{r(v-1)}{2}$, $3 \leq d \leq r$,  $r | dg$,  and 
$m=\frac{dg}{3}$, then there exists a $\KSS(v,m)$, provided 
$N \leq 3\left\lfloor \frac{g-1}{r-1} \right\rfloor$ where 
$N$ is defined by $\frac{rN}{3}=d(q+1+z)-(v-3m)$, 
$q=\left\lfloor \frac{v-3m}{3\left\lfloor 
\frac{d}{3}\right\rfloor} \right\rfloor$ and $0 \leq z < r$ 
is such that $3d(q+1+z) \equiv 0 \pmod{r}$.
\end{theorem}

\bigskip
The table below  illuminates existence results  
of Kirkman signal sets $\KSS(\texttt{v},\texttt{m})$ on
$\texttt{v}\equiv 1,3\pmod{6}$ points, $\texttt{v} \leq
3000$,
with partial parallel class size $\texttt{m}$.  
 The
parameter pair $\texttt{v},\texttt{m}$ is ommited when
$\texttt{m}$ divides $M$
and a $\KSS(\texttt{v},M)$ is known.
Paremeters pairs are anotated with ? if existence is unknown;
otherwise a letter 
\texttt{A},
\texttt{B},
\texttt{C},
\texttt{D},
\texttt{E},
\texttt{F}, or
\texttt{G},
is provided indicating which of the above theorems provides a constrution.

\verbatiminput{report3.txt}
\end{document}